**Michiel Hazewinkel**
Tel.: +31-35-6937033
fax: +31-35-6917401
homepage: michhaz.home.xs4all.nl
email: michhaz@xs4all.nl

# Left differential operators on non-commutative algebras

Michiel Hazewinkel
Burg. 's Jacob laan  18
NL-1401BR  BUSSUM
The Netherlands

1. **Introduction**. In [6 Grothendieck et al.], chapter 16, there is a comprehensive discussion of (higher order left) differential operators on schemes and hence on commutative algebras. Even more details can be found in [5 Gabriel], especially in connection with group schemes and their infinitesimal structure.

These discussions include a recursive description of these algebras of differential operators in terms of commutators with (left) multiplication operators; see propositions 16.8.8 and 16.8.9 on pages 42 and 43 of [6 Grothendieck et al.]. This description will be recalled later (sections 6, 7; for the affine case only).

More or less recently there has arisen substantial interest in whether something more or less similar can be done for non-commutative algebras. Straightforward generalizations definitely do not quite work.

The contexts in which this interest arises are, among more, Batalin-Vilkovisky algebras, Gerstenhaber algebras, homotopy algebras, higher derived algebras. See [1 Akman;  2 Akman;  3 Akman et al.] and, especially, the references in these papers. Also the deformation theory of algebras and diagrams à la Gerstenhaber-Shack, [4 Doubek]. The context 'homotopy algebras' is currently probably the most important. For some of this look at the last sections of [9 Markl].

There is in addition interest in such algebras in the realm of formal groups; for that consult the chapter on 'divided power algebras and sequences and algebras of differential operators' in [8 Hazewinkel].

The setting for this whole note is that of associative algebras over a base ring $k$ that is commutative and unital. Often it is good to think of $k$ as the ring of integers $\mathbf{Z}$.

2. **Example**. **Algebras of differential operators on $k[X,Y]$**. Certainly a very natural algebra of differential operators on $k[X,Y]$ is formed by the collection of operators that are finite sums of the form

$$\sum f_{i,j} \partial_X^i \partial_Y^j, \ f_{i,j} \in k[X,Y], \ i,j \geq 0 \tag{2.1}$$

Let me point out a few features of this algebra, which I will denote $\mathit{Diff}_{naive}(k[X,Y])$ [1]. First of all it is indeed a unital associative (and non-commutative) algebra over $k$. Second, there is a notion of degree. An expression like (2.1) is of degree $\leq n$ if $f_{i,j} = 0$ for all $i+j > n$. And this notion satisfies that if $D$ and $D'$ are of degree $\leq n, \leq n'$, respectively, then $DD'$ is of degree $\leq n+n'$. Third, it contains all derivations and obvious higher derivations on $k[X,Y]$. Fourth, if $D$ is an element of $\mathit{Diff}_{naive}(k[X,Y])$ of degree $\leq n$ and $t$ is an element of $k[X,Y]$, i.e. a polynomial, then the commutator

---
1. The reason for the subscript 'naive' will become clear later.

Michiel Hazewinkel
Tel.: +31-35-6937033
fax: +31-35-6917401
homepage: michhaz.home.xs4all.nl
email: michhaz@xs4all.nl

$$Dl_t - l_t D \tag{2.2}$$

is an element of $Diff_{naive}(k[X,Y])$ of degree $\leq n-1$. Here $l_t$ stands for the multiplication operator $l_t : f \mapsto tf$ (which nicely agrees with the notation used in (2.1))[2]. Fifth, the multiplications with an element from $k[X,Y]$ are in $Diff_{naive}(k[X,Y])$. This is natural in view of (2.2) albeit not strictly necessary.

These last three properties, with the fourth one playing a central role, permit one to define a nice generalization $Diff(A)$ of $Diff_{naive}(k[X,Y])$ for any algebra $A$.

There are of course other important algebras of differential operators on $k[X,Y]$, notably $k[\partial_X, \partial_Y]$, the algebra of differential operators with constant coefficients.[3]

3. **Infinitesimal neighborhoods**. **Infinitesimal thickenings**. The idea of infinitesimal neighborhoods (infinitesimal thickenings) is one guiding principle for both algebras of derivations and algebras of differentiable operators.

In this section the $k$-algebra $A$ is assumed commutative.

For the moment I will use a somewhat geometric language. So let $X = \text{Spec}(A)$ be the affine scheme (variety) (over $S = \text{Spec}(k)$) defined by the $k$-algebra $A$, and let it be the closed subscheme of an affine scheme $Y = \text{Spec}(B)$ defined by an ideal $J$ of $B$ (a prime ideal if $X$ is irreducible and reduced).

Then the infinitesimal neighborhoods of $X$ in $Y$ are the $X^{(n)} = \text{Spec}(B/J^{n+1})$. The underlying spaces of the $X^{(n)}$ are the same as that of $X$. But the structure sheaf is larger; the nilpotents $J/J^{n+1}$ have been added so to speak. Whence the phrases infinitesimal neighborhood and infinitesimal thickening.

In addition one defines $X^{(\infty)} = \text{Spec}(\varprojlim (B/J^{n+1}))$.

For instance if $x$ is a closed point in $Y$, defined by the maximal ideal $\mathfrak{m}_x$, the first infinitesimal thickening is the point $x$ together with the tangent space at $x$ to $Y$. Which is, as Shafarevich writes in [11 Shafarevich], section V.3.4, quite a large subscheme of $Y$.

There are two more or less canonical families of infinitesimal neighborhoods that play an important role in the matter of differential operators and derivations.

The first very easy one concerns $X = \text{Spec}(A)$ as a closed subscheme of the affine line over it; i.e the imbedding of schemes given by the algebra morphism $A[[t]] \longrightarrow A, \ t \mapsto 0$. then the infinitesimal neighborhoods are the $\text{Spec}(A[t]/(t^n))$ and $\text{Spec}(A[[t]])$.

---

2. The '$t$' stand for 'test'. The operator $D$ is *tested* by forming the commutator (2.2).
3. But these make little sense in situations for which the polynomials (or power series, or …) are a local model. Think of differential operators on (the smooth functions on) differentiable manifolds for example.

**Michiel Hazewinkel**
Tel.: +31-35-6937033
fax: +31-35-6917401
homepage: michhaz.home.xs4all.nl
email: michhaz@xs4all.nl

The second family consists of the infinitesimal neighborhoods of the diagonal $X \subset X \times X$. At the algebra level the diagonal imbedding is given by the multiplication morphism $m: A \otimes A \longrightarrow A$.[4] Let $J$ be the kernel of $m$ and define $P_A^n = (A \otimes A)/J^{n+1}$. These algebras play a crucial role in the definition of differential operators on $A$ as will be shown in a moment.

4. **Differential operators on a commutative algebra**. Again let $A$ be a commutative algebra over a (commutative unital) base ring $k$.

There are two natural algebra morphisms $A \longrightarrow A \otimes A$, viz $a \mapsto a \otimes 1$ and $a \mapsto 1 \otimes a$ and correspondingly two algebra morphisms $A \longrightarrow P_A^n$. Let $j_A^n$ be the second of these

$$j_A^n : A \longrightarrow P_A^n, \quad a \mapsto 1 \otimes a \mod J^{n+1} \tag{4.1}$$

By definition the (left) differential operators of order $\leq n$ are now the $k$-module morphisms of the form $\varphi j_A^n$ where $\varphi : P_A^n \longrightarrow A$ is a morphism of left $A$-modules. More generally, if $N$ and $M$ are two left $A$-modules a differential operator $M \longrightarrow N$ of order $\leq n$ is a $k$-module morphisms of the form

$$\varphi(j_A^n \otimes \mathrm{id}_M) : M \xrightarrow{\simeq} A \otimes_A M \xrightarrow{j_A^n \otimes \mathrm{id}_M} P_A^n \otimes_A M \xrightarrow{\varphi} N$$

where $\varphi$ is a morphisms of left $A$-modules and where in forming the tensor product $P_A^n \otimes_A M$ the right $A$-module structure on $P_A^n$ is that given by $j_A^n$ and the left $A$-module structure on $P_A^n \otimes_A M$ comes from that on $P_A^n$ induced by the morphism $A \longrightarrow A \otimes A, \ a \mapsto a \otimes 1$.

5. **Example**: **differential operators of order** $\leq n$ **on** $k[X]$. To see that this fits with the established and well understood notion of differential operators on rings of polynomials look at the simplest example $A = k[X]$ of polynomials over $k$ in a single indeterminate $X$. Then $A \otimes A = k[X, X']$, $J = (X' - X)$, and $P_A^n = k[X, dX]/((dX)^{n+1})$ where $dX$ is simply shorthand for $X' - X$. Thus, as a left $A = k[X]$ module

$$P_A^n = k[X] \oplus k[X]dX \oplus \cdots \oplus k[X](dX)^n$$

and an $A$-module morphism $\varphi : P_A^n \longrightarrow A$ is given by specifying $(n+1)$ polynomials $f_i \in k[X] = A, \ \varphi((dX)^i) = f_i$. Now

---
4. This uses commutativity. If $A$ is not commutative the multiplication map needs not be a morphism of algebras.

Michiel Hazewinkel
Tel.: +31-35-6937033
fax: +31-35-6917401
homepage: michhaz.home.xs4all.nl
email: michhaz@xs4all.nl

$$X'^m = (X + dX)^m = \sum_{i=0}^{m} \binom{m}{i} X^{m-i}(dX)^i \tag{5.1}$$

and so the composed morphism $\varphi j_A^n$ takes $X^m$ to

$$X^m \mapsto X'^m \mapsto f_0 X^m + f_1 \binom{m}{1} X^{m-1} + f_2 \binom{m}{2} X^{m-2} + \cdots + f_n \binom{m}{n} X^{m-n}$$

(where of course $\binom{m}{n} = 0$ for $n > m$ which is precisely equal to

$$(f_0 + f_1 \frac{d}{dX} + f_2 \frac{d^2}{dX^2} + \cdots + f_n \frac{d^n}{dX^n})(X^m)$$

The fact that this works out so precisely, including the divisibility properties over the integers of the operators involved, gives very strong support that this is precisely the right definition of differential operator.

6. **Characterization of differential operators on a commutative algebra**. Given a commutative algebra $A$ over $k$ let $l_x$ denote the $k$-morphism of multiplication by the element $x$.[5]

Now observe that in the example described above in section 5, if $D = f_0 + f_1 \frac{d}{dX} + f_2 \frac{d^2}{dX^2} + \cdots + f_n \frac{d^n}{dX^n}$ is a differential operator of order $\leq n$ on $k[X] = A$, then the commutator

$$Dl_x - l_x D = \mathrm{ad}(x)D \tag{6.1}$$

is a differential operator of order $\leq (n-1)$. This works in general.

6.2. **Theorem**. Let $D$ be a $k$-module endomorphism of the commutative $k$-algebra $A$. Then $D$ is a differential operator of order $\leq n$ if and only if for any $(n+1)$-tuple of elements $x_0, x_1, \cdots, x_n \in A$ the iterated commutator

$$\mathrm{ad}(x_0)\mathrm{ad}(x_1) \cdots \mathrm{ad}(x_n)D = 0 \tag{6.3}$$

---

5. the symbol $l_x$ really stands for multiplication on the left by $x$. Here, because of commutativity, there is no difference with $r_x$, multiplication on the right. But in the non-commutative case, to be discussed shortly, the difference between the two is far reaching and most important.

**Michiel Hazewinkel**
Tel.: +31-35-6937033
fax: +31-35-6917401
homepage: michhaz.home.xs4all.nl
email: michhaz@xs4all.nl

Moreover the differential operators on $A$ form an algebra under composition. More precisely, if $D$ is a differential operator of order $\leq n$, and $D'$ is a differential operator of order $\leq n'$ then the composed operator $DD'$ is a differential operator of order $\leq n + n'$.

6.4. The filtered algebra of differential operators thus defined is denoted $\mathit{Diff}(A)$ and called the algebra of (left) differential operators on $A$. It is not difficult to check that on $\mathbf{Q}[X,Y]$ this is precisely the naïve algebra of differential operators $\mathit{Diff}_{naive}\mathbf{Q}[X,Y]$ consisting of all finite sums $\sum f_{i,j} \partial_X^i \partial_Y^j$. The algebra of differential operators on $\mathbf{Z}[X,Y]$ is larger than $\mathit{Diff}_{naive}\mathbf{Z}[X,Y]$. It consists of all finite sums $\sum (i!)^{-1}(j!)^{-1} f_{i,j} \partial_X^i \partial_Y^j$, $f_{i,j} \in \mathbf{Z}[X,Y]$.

6.5. **Bibliographical notes**. The material treated in sections 2 - 6 above is merely the affine case of what is in [6 Grothendieck et al.], §16 in greater depth and in the greater generality of schemes (and without examples). In particular theorem 6.2. is the affine case of propositions 16.8.8, 16.8.9 on pages 42, 43. For still more detail, especially in connection with group schemes and their infinitesimal structure, see [5 Gabriel].

Most of the literature I have seen takes the criterion (6.3) as the definition of a (left) differential operator of order $\leq n$.

The main purpose of this note is to suggest a version of the description (6.3) for not necessarily commutative algebras which reduces to the one of theorem 6.2 in the commutative case and has a similar flavor.

7. **Recursive description of differential operators on a commutative algebra**. Theorem 6.2 implies the following recursive description of $\mathit{Diff}(A)$ (for commutative algebras).

$$\begin{aligned}
&\mathcal{D}_{-1}(A) = 0 \\
&\mathcal{D}_0(A) = \{d \in \mathrm{End}_k(A) : dl_t - l_t d \in \mathcal{D}_0(A) \text{ for all } t \in A\} \\
&\quad \vdots \\
&\mathcal{D}_{n+1}(A) = \{d \in \mathrm{End}_k(A) : dl_t - l_t d \in \mathcal{D}_n(A) \text{ for all } t \in A\} \\
&\quad \vdots \\
&\mathit{Diff}(A) = \bigcup_n \mathcal{D}_n(A)
\end{aligned} \qquad (7.1)$$

Note that $\mathcal{D}_0(A) = \{l_t : t \in A\} = A$, which because of commutativity is the same as the ring of right multiplications $r_t : a \mapsto at$, $t \in A$.

8. **Categorical direct sum in the category of associative algebras**. Now let $A$ be a not necessarily commutative algebra. The problem at hand is to find good generalizations of the notion of differential

**Michiel Hazewinkel**
Tel.: +31-35-6937033
fax: +31-35-6917401
homepage: michhaz.home.xs4all.nl
email: michhaz@xs4all.nl

operators for not necessarily commutative (but associative) algebras. A first thing to do is try to find a good generalization of the notion of an diagonal imbedding for 'non-commutative spaces'. For commutative spaces, e.g. in affine algebraic geometry the notion of a diagonal is clear. The product of two affine spaces $X = \mathrm{Spec}(A)$ and $Y = \mathrm{Spec}(B)$ is $X \times Y = \mathrm{Spec}(A \otimes B)$. I.e. the notion of a product of spaces corresponds on the algebraic side to the tensor product of commutative algebras. That fits of course with the fact that the tensor product $A \otimes_k B$ in the category of commutative algebras over $k$, together with the canonical imbeddings $j_A : A \longrightarrow A \otimes_k B$, $x \mapsto x \otimes 1_B$ and $j_B : B \longrightarrow A \otimes_k B$, $y \mapsto 1_A \otimes y$ is the categorical direct sum in the category of commutative algebras over $k$.

This is no longer true for not necessarily commutative algebras. Instead one must take the so-called free products, also called coproducts. These exist, see [10 Rowen], page 76, theorem 1.4.30.

One result is that the free product of the free algebras $k\langle X_1, \cdots, X_m \rangle$ and $k\langle Y_1, \cdots, Y_n \rangle$ is the free algebra $k\langle X_1, \cdots X_m, Y_1, \cdots, Y_n \rangle$ which is intuitively just right[6].

Now proceed as before. Let $A * B$ (or $A \amalg B$, which is better really but notationally ugly) denote the free product of the $k$-algebras $A$ and $B$ and let $j_A : A \longrightarrow A * B$ and $j_B : B \longrightarrow A * B$ be the corresponding canonical algebra morphisms (which need not be injective in general), then the defining property of $(A * B, j_A, j_B)$ is

> For every associative algebra $C$ together with algebra morphisms $\varphi_A : A \longrightarrow C$
> and $\varphi_B : B \longrightarrow C$ there is a unique morphisms of algebras $\psi : A * B \longrightarrow C$ (8.1)
> such that $\psi j_A = \varphi_A$, $\psi j_B = \varphi_B$

For the case $A = B$ write $j_1 : A \mapsto A * A$ for "$j_A$" and $j_2 : A \longrightarrow A * A$ for "$j_B$". The codiagonal of an algebra $A$ is now the unique morphisms of algebras $\delta_{codiag} : A * A \longrightarrow A$ such that $\delta_{codiag} j_1 = \mathrm{id}_A$, $\delta_{codiag} j_2 = \mathrm{id}_A$.

In case $A = k\langle X_1, \cdots, X_m \rangle$, $A * A = k\langle X_1, \cdots, X_m; X'_1, \cdots, X'_m \rangle$ (using a somewhat obvious notation), $j_1(X_i) = X_i$, $j_2(X_i) = X'_i$ and the codiagonal is given by $\delta_{codiag}(X_i) = X_i = \delta_{codiag}(X'_i)$ with as kernel the two-sided ideal generated by the $X'_i - X_i$.

9. **Algebras over a non-commutative algebra**. Let $A$ and $B$ be algebras over the commutative base ring $k$. In these jottings an $A$-algebra $B$ is simply a $k$-algebra together with a morphism of $k$-algebras $\rho : A \longrightarrow B$. So the left and right action of an element $a$ of $A$ on an element of $b$ of $B$ are given by $ab = \rho(a)b$, $ba = b\rho(a)$.

---

6. This is all hindsight. I first messed around for more than a full day with tensor products of non-commutative algebras getting absolutely nowhere.

**Michiel Hazewinkel**
Tel.: +31-35-6937033
fax: +31-35-6917401
homepage: michhaz.home.xs4all.nl
email: michhaz@xs4all.nl

In much of the literature one finds the definition that an $A$-algebra is a ring $C$ which is also a left $A$-module such that $a(c_1 c_2) = (ac_1)c_2 = c_1(ac_2)$. That is different and practically forces the $A$ to be commutative in that such an algebra can as well be considered to be an algebra over the commutative ring $A_{ab} = A/J$ where $J$ is the ideal generated by the commutator differences $a_1 a_2 - a_2 a_1$.

The definition given above in the first paragraph of this section can be easily extended to cover not necessarily unital (but associative, at least for the current purposes) algebras. Such an algebra is an $A$-bimodule (i.e left and right multiplications of elements of $B$ by elements of $A$ are defined and $(a_1 b)a_2 = a_1(ba_2)$), and there is an $A$-bilinear multiplication $(b_1, b_2) \mapsto b_1 b_2$ such that $a_1 b_1 a_2 b_2 \cdots a_n b_n a_{n+1}$ is one and the same element of $B$ whatever bracketing is used to calculate it.

10. **Multimorphisms**. Let $B$ and $C$ be $A$-algebras (in the sense above, not necessarily unital). A $k$-morphism $\varphi : B \mapsto C$ is an $A$-multi-morphism if it satisfies

$$\varphi(a_1 b_1 a_2 b_2 \cdots a_n b_n a_{n+1}) = a_1 \varphi(b_1) a_2 \varphi(b_2) \cdots a_n \varphi(b_n) a_{n+1} \qquad (10.1)$$

for all $a_1, a_2, \cdots a_{n+1} \in A$, $b_1, b_2, \cdots, b_n \in B$

Obviously any $A$-algebra morphism from $B$ to $C$ is an $A$-multimorphism, but the notion is more general than that.

11. **Example**. Consider the free algebra $B = k\langle X; Y \rangle = k\langle X_1, \cdots, X_m; Y_1, \cdots, Y_n \rangle$ in the noncommuting indeterminates $X_1, \cdots, X_m; Y_1, \cdots, Y_n$ and consider it as algebra over $A = k\langle X \rangle$ via the obvious natural embedding[7]. Now let $C$ be any $A$-algebra. For every $Y_i$ let $z_i$ be an element from $C$ and let $B^{(r)}$ be the (non-unital) sub-$A$-module of $B$ spanned by all monomials which have precisely $r$ $Y$'s in them. Every such monomial can be written uniquely in the form

$$M_1 Y_{i_1} M_2 Y_{i_2} \cdots M_r Y_{i_r} M_{r+1}$$

where the $M_i$ are monomials in the $X$'s (including the "empty" monomial 1). Now define

$$\varphi(M_1 Y_{i_1} M_2 Y_{i_2} \cdots M_r Y_{i_r} M_{r+1}) = M_1 z_{i_1} M_2 z_{i_2} \cdots M_r z_{i_r} M_{r+1}$$

on these monomials and $\varphi(B^{(s)}) = 0$ for $s \neq r$.

Then $\varphi : B' \longrightarrow C$ is an $A$-multimorphism (but not an $A$-algebra morphism).

---

7. Nota Bene: $k\langle X; Y \rangle$ is not the same as $k\langle X \rangle\langle Y \rangle$; the latter is obtained from the former by quotienting out the ideal generated by the commutators $X_i Y_j - Y_j X_i$. $i = 1, \cdots, m$, $j = 1, \cdots, n$.

**Michiel Hazewinkel**
Tel.: +31-35-6937033
fax: +31-35-6917401
homepage: michhaz.home.xs4all.nl
email: michhaz@xs4all.nl

12. **Differential operators on a non-commutative algebra.** A first potential definition of a differential operators on an associative, but not necessarily commutative, algebra of order $\leq n$ now proceeds as follows. For each $a \in A$ let $l_a \in \mathrm{End}_k(A)$ denote the $k$-module endomorphism of multiplication on the left $l_a(x) = ax$; multiplications on the right will also be needed: $r_a(x) = xa$. The subalgebra of $\mathrm{End}_k(A)$ of left multiplication morphisms will be denoted $L_A$; that of right multiplication morphisms $R_A$.

It seems reasonable to require at least two properties of the to be defined algebra of differential operators $\mathcal{D}(A)$. First that it contain all the terms that occur in the Hasse-Schmidt derivations on $A$ and second that it be a multiplicatively filtered algebra.

As in the commutative case set

$$\mathcal{D}_{-1}(A) = 0 \tag{12.1}$$

Next ( proceeding in analogy with the commutative case consider commutators with left multiplication operators) define

$$\mathcal{D}'_0(A) = \{\varphi \in \mathrm{End}_k(A) : l_t\varphi - \varphi l_t \in \mathcal{D}_{-1}(A) = 0\} \tag{12.2}$$

Then for $\varphi \in \mathcal{D}'_0$ for all $t \in A$,

$$t\varphi(1) = (l_t\varphi)(1) = (\varphi l_t)(1) = \varphi(t)$$

so that $\varphi$ is $r_{\varphi(1)}$, right multiplication by $\varphi(1)$. Inversely all right multiplication morphisms satisfy the condition in (12.2). So $\mathcal{D}'_0 = R_A$.

Now lets take a look at the requirement that the module morphisms occurring in a Hasse-Schmidt derivation occur in the to be defined algebra. For material on Hasse-Schmidt derivations see [7 Hazewinkel] and the references therein. One has for a higher derivation $\partial_n$ occurring in a Hasse-Schmidt derivation sequence $(\partial_0 = \mathrm{id}, \partial_1, \partial_2, \cdots)$

$$\partial_n(tx) = t\partial_n(x) + \partial_1(t)\partial_{n-1}(x) + \cdots + \partial_{n-1}(t)\partial_1(x) + \partial_n(t)x$$

so that

$$\partial_n l_t - l_t \partial_n = l_{\partial_1(t)}\partial_{n-1} + \cdots + l_{\partial_{n-1}(t)}\partial_1 + l_{\partial_n(t)} \tag{12.3}$$

Intuitively at least the to be defined $\mathcal{D}_n(A)$ should contain the $\partial_n$ of a Hasse-Schmidt derivation. So formula (12.3) strongly suggests that the $\mathcal{D}_n(A)$ be stable under composition with a left multiplication operator. Looking in particular at $n = 1$ in equation (12.3) it seems inevitable that $\mathcal{D}_0(A)$ be equal to

**Michiel Hazewinkel**
Tel.: +31-35-6937033
fax: +31-35-6917401
homepage: michhaz.home.xs4all.nl
email: michhaz@xs4all.nl

$$\mathcal{D}_0(A) = L_A \mathcal{D}'_0(A) = L_A R_A = R_A L_A$$

(where the last equation holds because left multiplications commute with right multiplications), and where $L_A R_A$ is to be interpreted as all finite sums of expressions of the form $l_a r_b$. Now it is also desired that the $\mathcal{D}_n(A)$ define a multiplicative filtration. That implies in particular that $\mathcal{D}_n(A)\mathcal{D}_0(A) \subset \mathcal{D}_n(A)$, so that the $\mathcal{D}_n(A)$ also need to be stable under composition on the right with a left multiplication operator.

Thus one is led to the following inductive definition

$$\begin{aligned}
\mathcal{D}_{-1}(A) &= 0 \\
\mathcal{D}'_0(A) &= \{\varphi \in \mathrm{End}_k(A) : l_t\varphi - \varphi l_t \in \mathcal{D}_{-1}(A) = 0\} = R_A \\
\mathcal{D}_0(A) &= R_A L_A = L_A R_A L_A = L_A \mathcal{D}'_0(A) L_A \\
\mathcal{D}'_1(A) &= \{\varphi \in \mathrm{End}_k(A) : l_t\varphi - \varphi l_t \in \mathcal{D}_0(A)\} \\
\mathcal{D}_1(A) &= L_A \mathcal{D}'_1(A) L_A \\
&\vdots \\
\mathcal{D}'_n(A) &= \{\varphi \in \mathrm{End}_k(A) : l_t\varphi - \varphi l_t \in \mathcal{D}_{n-1}(A)\} \\
\mathcal{D}_n(A) &= L_A \mathcal{D}'_n(A) L_A \\
&\vdots
\end{aligned} \qquad (12.4)$$

where for two subrings $\mathcal{D}$ and $\mathcal{D}'$ the expression $\mathcal{D}\mathcal{D}'$ is to be interpreted as all finite sums of products $D_1 D_2$, $D_1 \in \mathcal{D}, D_2 \in \mathcal{D}'$.

The operators in $\mathcal{D}_n(A)$ are called (left) differential operators of order $\leq n$ and one easily checks that $\mathcal{D}_{n-1}(A) \subset \mathcal{D}_n(A)$. The full ring of differential operators on $A$ is now defined as the inductive limit of the $\mathcal{D}_n(A)$

$$\mathcal{D}(A) = \bigcup_n \mathcal{D}_n(A)$$

The union being taken in $\mathrm{End}_k(A)$.

13. **Theorem**. The filtration of $\mathcal{D}(A)$ by the $\mathcal{D}_n(A)$ is multiplicative; i.e. $\mathcal{D}_r(A)\mathcal{D}_s(A) \subset \mathcal{D}_{r+s}(A)$.

Proof. This is done by induction (and rather tedious calculations). The start of the induction is the proof that

$$\mathcal{D}_0(A)\mathcal{D}_n(A) \subset \mathcal{D}_n(A) \quad \text{and} \quad \mathcal{D}_n(A)\mathcal{D}_0(A) \subset \mathcal{D}_n(A) \qquad (13.1)$$

Michiel Hazewinkel
Tel.: +31-35-6937033
fax: +31-35-6917401
homepage: michhaz.home.xs4all.nl
email: michhaz@xs4all.nl

I shall only write out the details for the first of these inclusions (equalities in fact); the second one goes in exactly the same way. Because $\mathfrak{D}_n(A)$ is stable under composition with left multiplication operators on the left and on the right it suffices to prove that

$$r_a D \in \mathfrak{D}_n(A) \text{ for } D \in \mathfrak{D}'_n(A)$$

This is clear for $n = 0$. So induction can be used. Now

$$l_t r_a D - r_a D l_t = r_a l_t D - r_a D l_t \in r_a \mathfrak{D}_{n-1}(A) \subset \mathfrak{D}_{n-1}(A)$$

where the last inclusion is by the induction hypothesis. This proves the first inclusion of (13.1).

To finish the proof one shows that

$$\mathfrak{D}_r(A)\mathfrak{D}_s(A) \subset \mathfrak{D}_n(A), \text{ where } n = r + s$$

Because $\mathfrak{D}_n(A)$ is stable under composition with left multiplications on the left and on the right and because $l_a l_b = l_{ab}$ it suffices to show that

$$D_1 l_a D_2 \in \mathfrak{D}_n(A) \text{ in the situation } D_1 \in \mathfrak{D}'_r(A), D_2 \in \mathfrak{D}'_s(A) \tag{13.2}$$

That means examining the expression

$$D_1 l_a D_2 l_t - l_t D_1 l_a D_2$$

and showing that it is in $\mathfrak{D}_{n-1}(A)$. The procedure is to use commutators and induction to move the left multiplication operators to the left. Here are the (tedious) details.

$$D_1 l_a D_2 l_t - l_t D_1 l_a D_2 = D_1 l_a D_2 l_t - D_1 l_a l_t D_2 + D_1 l_a l_t D_2 - l_t D_1 l_a D_2 \tag{13.3}$$

Now $D_2 l_t - l_t D_2 \in \mathfrak{D}_{r-1}(A)$ and so with induction the first two summands of (13.3) together form an element of $\mathfrak{D}_{n-1}(A)$. Next

$$D_1 l_a l_t D_2 - l_t D_1 l_a D_2 = D_1 l_a l_t D_2 - D_1 l_t l_a D_2 + D_1 l_t l_a D_2 - l_t D_1 l_a D_2 \tag{13.4}$$

This time the last two terms combine to form an element of $\mathfrak{D}_{n-1}(A)$ because $D_1 l_t - l_t D_1 \in \mathfrak{D}_{r-1}(A)$. Next

$$D_1 l_a l_t D_2 - D_1 l_t l_a D_2 = D_1 l_a l_t D_2 - l_a D_1 l_t D_2 + l_a D_1 l_t D_2 - D_1 l_t l_a D_2$$

and the first two terms combine to form an element of $\mathfrak{D}_{n-1}(A)$. Next

$$l_a D_1 l_t D_2 - D_1 l_t l_a D_2 = l_a D_1 l_t D_2 - l_a l_t D_1 D_2 + l_a l_t D_1 D_2 - D_1 l_t l_a D_2$$

Michiel Hazewinkel
Tel.: +31-35-6937033
fax: +31-35-6917401
homepage: michhaz.home.xs4all.nl
email: michhaz@xs4all.nl

and because $D_1 l_t - l_t D_1$ is in $\mathcal{D}_{r-1}(A)$, the first two term combine to form an element of $\mathcal{D}_{n-1}(A)$. Next

$$l_a l_t D_1 D_2 - D_1 l_t l_a D_2 = l_a l_t D_1 D_2 - l_t D_1 l_a D_2 + l_t D_1 l_a D_2 - D_1 l_t l_a D_2$$

This time the last two summands combine to form an element of $\mathcal{D}_{n-1}(A)$. Next

$$l_a l_t D_1 D_2 - l_t D_1 l_a D_2 = l_a l_t D_1 D_2 - l_t l_a D_1 D_2 + l_t l_a D_1 D_2 - l_t D_1 l_a D_2$$

Again the last two summands combine to give an element of $\mathcal{D}_{n-1}(A)$. Finally

$$l_a l_t D_1 D_2 - l_t l_a D_1 D_2 = l_{[a,t]} D_1 D_2$$

and so, again because $\mathcal{D}_n(A)$ is stable under composition with left multiplication operators on the left, things are reduced to proving that

$$D_1 D_2 \in \mathcal{D}_n(A), \text{ when } D_1 \in \mathcal{D}'_r(A), \ D_2 \in \mathcal{D}'_s(A), \ n = r+s$$

Now

$$D_1 D_2 l_t - l_t D_1 D_2 = D_1 D_2 l_t - D_1 l_t D_2 + D_1 l_t D_2 - l_t D_1 D_2 \tag{13.5}$$

Both the first and second summands on the one hand and the third and fourth summands on the other combine to give an element of $\mathcal{D}_{n-1}(A)$ (again by induction).

This concludes the proof.

14. **Concluding remarks**. There remain a number of things to be sorted out in the future. Notably:

1) What happens if 'left' is replace by 'right' (and vice versa) everywhere.

2) To work out by way of example what things look like in the case that $A$ is a free algebra. For higher derivations (Hasse-Schmidt derivations), a related matter by section 12 above, this has been done in the chapter on divided power algebras in [8 Hazewinkel]. Also there are there constructions which look very much like multimorphisms and they play a significant role.

3) To obtain a description of $Diff(A)$ as described in section 12 in terms of infinitesimal thickenings and the like. I think that free products and multimorphisms can play a significant role in that, whence the inclusion of sections 8 and 10 in the above.

**References**

1. Akman, Füsun, *On some generalizations of Batalin-Vilkoovisky algebras*, 1996, arXiv:q-alg/9506027 v3.

Michiel Hazewinkel
Tel.: +31-35-6937033
fax: +31-35-6917401
homepage: michhaz.home.xs4all.nl
email: michhaz@xs4all.nl

2.      Akman, Füsun, *Chicken or egg? A hierarchy of homotopy algebras*, Homology, homotopy and applications **7** (2005), 5 - 39.
3.      Akman, Füsun and Lucian M Ionescu, *Higher derived brackets and deformation theory I*, A 2008, arXiv:math/0504541 v2 [math.QA].
4.      Doubek, Martin, *On resolutions of diagrams of algebras*, 2011, arXiv:1107.1408 v1 [math.AT].
5.      Gabriel, P, *Étude infinitésimale des schémas en groupes. A: Opérateurs différentielles et p-algèbres de Lie*, in M. Demazure and A. Grothendieck, eds., *Séminaire de géométrie algébrique du Bois Marie 1962/1964, SGA3: Schémas en groupes I*, Springer, 1970, Exposé VIIA, 409 - 473.
6.      Grothendieck, Alexandre and Jean Dieudonné, *Éléments de géométrie algébrique IV. Étude local des schémas et morphismes de schémas, quatrième partie (EGA IV, 4)*, 1967.
7.      Hazewinkel, Michiel, *Hasse-Schmidt derivations and the Hopf algebra of non-commutative symmetric functions*, Axioms **2012** (2012), 149 - 154. Preprint version: arXiv:1110.6108.
8.      Hazewinkel, Michiel, *Formal groups and applications 2*, In preparation.
9.      Markl, Martin, *Loop homotopy algebras in closed string field theory*, 1999. arXiv:hep-th/9711045 v2.
10.     Rowen, Louis H, *Ring theory. Volume I*, Academic Press, 1988.
11.     Shafarevich, I R, *Basic algebraic geometry*, Springer, 1994.